\documentclass[a4paper,12pt]{amsart}
\usepackage{amssymb}
\usepackage{ifthen}
\nonstopmode \numberwithin{equation}{section}
\newtheorem{thm}[equation]{Theorem}
\newtheorem{cor}[equation]{Corollary}
\newtheorem{lem}[equation]{Lemma}
\newtheorem{prop}[equation]{Proposition}

\newtheorem{conj}{Conjecture}

\theoremstyle{definition}
\newtheorem{defn}{Definition}[section]

\newtheorem{prob}[equation]{Problem}

\newenvironment{rem}{%
\bigskip
\noindent \textsl{{\sl Remark. }}}{\bigskip}
\newenvironment{rems}{%
\bigskip
\noindent \textsl{{\sl Remarks. }}}{\bigskip}

\newcounter {own}
\def\theown {\thesection       .\arabic{own}}

\newenvironment{pf}[1][]{%
 \vskip 3mm
 \noindent
 \ifthenelse{\equal{#1}{}}%
  {{\bf Proof. }}%
  {{\bf #1.} }%
 }%
{\qed\bigskip}

\newcounter{alphabet}
\newcounter{tmp}

\newcommand{\IR}{{\mathbb R}}
\newcommand{\ID}{{\mathbb D}}
\newcommand{\IN}{{\mathbb N}}
\newcommand{\IC}{{\mathbb C}}

\def\be{\begin{equation}}
\def\ee{\end{equation}}

\newcommand{\bee}{\begin{enumerate}}
\newcommand{\eee}{\end{enumerate}}

\newcommand{\pay}{\!\!\!}
\newcommand{\blem}{\begin{lem}}
\newcommand{\elem}{\end{lem}}
\newcommand{\bthm}{\begin{thm}}
\newcommand{\ethm}{\end{thm}}
\newcommand{\bcor}{\begin{cor}}
\newcommand{\ecor}{\end{cor}}
\newcommand{\beg}{\begin{examp}}
\newcommand{\eeg}{\end{examp}}
\newcommand{\begs}{\begin{examples}}
\newcommand{\eegs}{\end{examples}}
\newcommand{\bdefe}{\begin{defin}}
\newcommand{\edefe}{\end{defin}}
\newcommand{\bprob}{\begin{prob}}
\newcommand{\eprob}{\end{prob}}
\newcommand{\bei}{\begin{itemize}}
\newcommand{\eei}{\end{itemize}}

\newcommand{\bcon}{\begin{conj}}
\newcommand{\econ}{\end{conj}}
\newcommand{\bcons}{\begin{conjs}}
\newcommand{\econs}{\end{conjs}}
\newcommand{\bprop}{\begin{prop}}
\newcommand{\eprop}{\end{prop}}
\newcommand{\br}{\begin{rem}}
\newcommand{\er}{\end{rem}}
\newcommand{\brs}{\begin{rems}}
\newcommand{\ers}{\end{rems}}
\newcommand{\bo}{\begin{obser}}
\newcommand{\eo}{\end{obser}}
\newcommand{\bos}{\begin{obsers}}
\newcommand{\eos}{\end{obsers}}
\newcommand{\bpf}{\begin{pf}}
\newcommand{\epf}{\end{pf}}
\newcommand{\ba}{\begin{array}}
\newcommand{\ea}{\end{array}}
\newcommand{\beq}{\begin{eqnarray}}
\newcommand{\beqq}{\begin{eqnarray*}}
\newcommand{\eeq}{\end{eqnarray}}
\newcommand{\eeqq}{\end{eqnarray*}}

\newcommand{\ra}{\rightarrow}

\newcommand{\ds}{\displaystyle}

\begin{document}
\bibliographystyle{amsplain}
\title{Characterization and the pre-Schwarzian norm estimate for concave univalent functions}
\author{B. Bhowmik}
\address{B. Bhowmik, Department of Mathematics,
Indian Institute of Technology Madras, Chennai-600 036, India.}
\email{ditya@iitm.ac.in}
\author{S. Ponnusamy}
\address{S. Ponnusamy, Department of Mathematics,
Indian Institute of Technology Madras, Chennai-600 036, India.}
\email{samy@iitm.ac.in}
\author{K.-J. Wirths}
\address{K.-J. Wirths, Institut f\"ur Analysis, TU Braunschweig,
38106 Braunschweig, Germany}
\email{kjwirths@tu-bs.de}

\date{ %\today
Version: Jan. 06, 08; File: bp7a-new${}_{-}$final.tex}
\begin{abstract}
Let $Co(\alpha)$ denote the class of concave univalent functions
in the unit disk  $\ID$. Each function  $f\in Co(\alpha)$ maps the
unit disk  $\ID$ onto the complement of an unbounded convex set. In this paper we find
the exact disk of variability for the functional $(1-|z|^2)\left ( f''(z)/f'(z)\right)$,
$f\in Co(\alpha)$.
In particular, this gives  sharp upper and lower estimates for
the pre-Schwarzian norm of concave univalent functions.
Next we obtain the set of
variability of the functional $(1-|z|^2)\left(f''(z)/f'(z)\right)$, $f\in Co(\alpha)$
whenever $f''(0)$ is fixed.
We also give a characterization for concave functions in terms of Hadamard convolution.
In addition to sharp coefficient inequalities, we prove that functions in
$Co(\alpha)$ belong to the $H^p$ space for $p<1/\alpha$.
%and obtain a test
%for concavity of analytic functions in terms of their Taylor coefficients.
\end{abstract}

\thanks{}

\maketitle
\pagestyle{myheadings}
\markboth{B. Bhowmik, S. Ponnusamy and K.-J. Wirths}{Norm estimates}

\section{Introduction and Preliminary Results}
Let $\mathcal{H}$ denote the class of functions analytic in the unit disk
$\ID := \{ z\in \IC: \,  |z| < 1 \}$. We denote the
class of locally univalent functions by  $\mathcal{LU}$. The class of locally
univalent functions is a vector space with respect to Hornich operations
(see \cite{H}).
For $f\in  \mathcal{LU}$, the pre-Schwarzian derivative $T_f$ is defined by
$$
T_f= \frac{f''}{f'}
$$
and we define the norm of $T_f$ by
$$
\|T_f\|= \sup_{z\in \ID}(1-|z|^2)|T_f(z)|.
$$
This is indeed a norm with respect to Hornich operations. It is known that
$\|f\|<\infty$ if and only if $f$ is uniformly locally univalent, i.e. there exists a
constant $r=r(f)>0$ such that $f$ is univalent in each disk of hyperbolic radius $r$ in
$\ID.$

Let $\mathcal{A}$ denote the class of functions $f\in {\mathcal H}$ with the normalization
$f(0)=f'(0)-1=0$ and $\mathcal{S}$ be the class of functions in $\mathcal{A}$ that are
univalent in $\ID$ and ${\mathcal H}_1$ denotes the class of functions $f$ in
$\mathcal H$ such that $f(0)=1$. Also we define the subclass
$\mathcal{K}\subset\mathcal{S}$ of
convex functions whenever $f(\ID)$ is a convex domain and the subclass $\mathcal{S}^*$
of starlike functions whenever $f(\ID)$ is a domain that is starlike with respect
to the origin (cf. \cite{Duren,P}).
It is well known that
$\|T_f\|\leq 6$ for $f \in \mathcal{S}$, and $\|T_f\|\leq 4$ for $f \in \mathcal{K}$.
Conversely, by Becker's theorem (\cite{Bek-72})
it follows that if  $f\in\mathcal{A}$ and $\|T_f\|\leq 1$ then $f \in \mathcal{S}$.

A function $f:\ID\to \IC$ is said to belong to the family $Co(\alpha)$  if
$f$ satisfies  the following conditions:
\bee
\item[(i)] $f$ is analytic in $\ID$ with the standard normalization $f(0)=f'(0)-1=0$.
In addition it satisfies $f(1)=\infty$.
\item[(ii)] $f$ maps $\ID$ conformally onto a set whose complement with respect to $\IC$ is
convex.
\item[(iii)] the opening angle of $f(\ID)$ at $\infty$ is less than or equal to
$\pi\alpha$, $\alpha\in (1,2]$.
\eee
This paper concerns the family $Co(\alpha)$ and in order to proceed with our investigation,
we recall the analytic characterization for functions in $Co(\alpha), \alpha \in (1,2]$:
$f\in Co(\alpha)$ if and only if
\be\label{p8eq1}
{\rm\, Re}\, P_f(z)> 0, ~z\in \ID,
\ee
where
%\be\label{p8eq1}
$$
P_f(z)=\frac{2}{\alpha-1}\,
 \left [\frac{(\alpha+1)}{2}\frac{1+z}{1-z}-1-z \frac{f''(z)}{f'(z)}\right ].
$$
The class $Co(\alpha)$ is referred to as the class of concave univalent functions
and for a detailed discussion about concave functions, we refer to
\cite{Avk-Wir-06, Avk-Wir-05,Pom-Cruz}.  We note
that for  $f\in Co(\alpha)$, $\alpha\in (1,2]$, the closed set $\IC\backslash f(\ID)$
is convex and unbounded. We observe that $Co(2)$ contains the classes $Co(\alpha)$,
$\alpha\in (1,2]$.

In this paper, we first find the exact set of variability
for the functional $(1-|z|^2)T_f(z)$ and as a consequence of this we derive
upper and lower bounds for the pre-Schwarzian norm
$\|T_f\|$, for functions $f$ in $Co(\alpha)$. Next we obtain the set of
variability of the functional $(1-|z|^2)T_f(z), f\in Co(\alpha)$
whenever $f''(0)$ is fixed.
Also,  we give a representation formula in terms of Hadamard
convolution for functions in $Co(\alpha)$ and some interesting link with the Kaplan
class. Lastly, we present  sharp inequalities among coefficients of functions in
$Co(\alpha)$.
%and using this we give a criteria for analytic functions to belong in $Co(\alpha)$ in
%terms of their Taylor coefficients.

\section{ Main results}\label{p8sec1}
First we prove the following lemma:
\blem\label{p7lem1}
Let $\psi \in \mathcal{H}_1$ be such that it is starlike with respect to $1$ and
suppose that $g\in \mathcal{A}$ satisfies
$$
\frac{2}{\alpha-1}\, \left [\frac{(\alpha+1)}{2}\frac{1+z}{1-z}-1-z \frac{g''(z)}
{g'(z)}\right ] = \psi(z), \quad z\in \ID,
$$
for some $\alpha\in (1,2]$.
Then, for $f\in Co(\alpha)$, the condition
\be\label{p7eq5}
\frac{2}{\alpha-1}\, \left [\frac{(\alpha+1)}{2}\frac{1+z}{1-z}-1-z \frac{f''(z)}
{f'(z)}\right ] \prec \psi(z)
\ee
implies $(1-z)^{\alpha + 1} f'(z) \prec (1-z)^{\alpha+1} g'(z)$.
\elem
\bpf
We first note that $\psi(0)=1$ and
$$
\psi'(0)=\frac{2}{\alpha-1}[(\alpha+1)-g''(0)],
$$
as $\psi$ is starlike and hence univalent. Also, we note that
$$
g'(z)=\exp\int_{0}^{z}\frac{(\alpha-1)[1-(1-\zeta)\psi(\zeta)]+(\alpha+3)\zeta}
{2\zeta(1-\zeta)}\,d\zeta, ~ z\in \ID,
$$
which is a non-vanishing analytic function in the unit disk.
Let
\beqq
h(z) &=& \frac{2}{(\alpha-1) \psi'(0)}\left[ -\log ((1-z)^{\alpha +1}g'(z))\right]\\
&& := - c \log ((1-z)^{\alpha +1}g'(z)), \quad c= \frac{2}{(\alpha-1) \psi'(0)}.
\eeqq
Since $\frac{(\alpha-1)c}{2} (\psi -1)\in \mathcal{A}$ is starlike, a computation
shows that
$$
1 + z \frac{h''(z)}{h'(z)}= \frac{z\psi'(z)}{\psi(z)-1}
$$
has positive real part and so $h(z)$ is convex with $h(0)=0=h'(0)-1$.
The condition (\ref{p7eq5}) and a little computation  reveals that
$$
c\left [(\alpha+1)\frac{z}{1-z}-z \frac{f''(z)}
{f'(z)}\right ] \prec c\left [(\alpha+1)\frac{z}{1-z}-z \frac{g''(z)}
{g'(z)}\right ] = z h'(z).
$$
Equivalently the above can be written as
$$
z\left[-c \log\left((1-z)^{\alpha +1}f'(z)\right)\right]'\prec z h'(z).
$$
As $h(z)$ is convex, by using a result due to
Suffridge \cite[p.~76, Theorem 3.1d]{Mil-Mocanu}, we get
$$
-c \log\left((1-z)^{\alpha +1}f'(z)\right)\prec h(z)=
-c \log\left((1-z)^{\alpha +1}g'(z)\right),
$$
which gives the desired result.
\epf

We now recall that, for $f,g\in \mathcal{A}$, the condition $f'\prec g'$
implies the inequality $\|T_f\|\leq \|T_g\|$ (see (\cite{Kim-2002})).
Hence we obtain

\bthm
Let $g$ be as Lemma {\rm \ref{p7lem1}}. If $f\in Co(\alpha)$, then
$\|T_F\|\leq \|T_G\|$ where
$$
F(z)=\int_{0}^{z}(1-\zeta)^{\alpha+1}f'(\zeta)\,d\zeta ~\mbox{and} ~
G(z)=\int_{0}^{z}(1-\zeta)^{\alpha+1}g'(\zeta)\,d\zeta.
$$
\ethm

Now we state the following corollary:

\bcor
For $f\in Co(\alpha)$ and $g$  as in Lemma {\rm \ref{p7lem1}}, we have
$$\left|(1-|z|^2)\frac{f''(z)}{f'(z)}-(\alpha+1)\frac{1-|z|^2}{1-z}\right|\leq
\left|(1-|z|^2)\frac{g''(w(z))}{g'(w(z))}-(\alpha+1)\frac{1-|w(z)|^2}{1-w(z)}\right|,
$$
where $w : \ID \ra \overline{\ID}$ is a holomorphic function with $w(0)=0$.
Equality holds when $w(z)=z$.
\ecor
\bpf
From Lemma \ref{p7lem1}, we have
$$
(1-z)^{\alpha + 1} f'(z) \prec (1-z)^{\alpha+1} g'(z).
$$
Using the definition of subordination we have,
$$
f'(z)=\frac{(1-w(z))^{\alpha+1}g'(w(z))}{(1-z)^{\alpha+1}},
$$
where $w : \ID \ra \overline{\ID}$ is a holomorphic function with $w(0)=0$.
After taking the logarithmic derivative and using the following
Schwarz-Pick inequality
$$|w'(z)|\leq \frac{1-|w(z)|^2}{1-|z|^2},
$$
we get the desired inequality stated in the corollary. Also it is
easy to see that equality holds in the inequality when $w(z)=z$.
\epf

Now, for  $f\in Co(\alpha),$  we  find the exact set of
variability for  the functional $(1-|z|^2)T_f(z)$,
which essentially gives both sharp upper and lower bounds for
the pre-Schwarzian norm $\|T_f\|$.

\bthm\label{p8thm1}
Let $\alpha \in (1,2]$ be fixed.
Then the set of variability of the functional $(1-|z|^2)T_f(z)$, $f\in Co(\alpha),$
is the closed disk with center
$$2\overline{z}\,+\,(\alpha + 1)(1-\overline{z})/(1-z)
$$
and radius $\alpha - 1$. The points on the boundary of this disk are attained if and only if
$f$ is one of the functions $g_{\theta}$,  where,
$$
g_{\theta}(z)=\frac{1}{\alpha(1+e^{i\theta})}
\left [\left(\frac{1+e^{i\theta}z}{1-z}\right)^\alpha-1 \right],\quad \mbox{for}~\,
\theta\in[0,2\pi]\setminus \{\pi\},
$$
and
$$
g_{\pi}(z)=\frac{z}{1-z}, \quad \mbox{for}~\, \theta=\pi.
$$
\ethm
\bpf
We use the characterization (\ref{p8eq1}) for functions in $Co(\alpha)$ and
the representation
$$P_f(z)\,=\,\frac{1\,-\,z\omega (z)}{1\,+\,z\omega (z)},
$$
where $\omega :\ID\to \overline{\ID}$ is an unimodular bounded analytic function.
It follows that
$$T_f(z)\,=\,\frac{(\alpha-1)\omega(z)+(\alpha+1)+2z\omega(z)}{(1-z)(1+z\omega(z))}.
$$
By a routine computation one recognizes that
$$ (1-|z|^2)T_f(z)\,-\,2\overline{z}\,-\,(\alpha + 1)
\frac{1-\overline{z}}{1-z}=(\alpha-1)\frac{\overline z+ \omega(z)}{1+z\omega(z)}.
$$
Hence, the condition $|\omega(z)|\leq 1$
is equivalent to
\be\label{p8eq5aa}
\left|(1-|z|^2)T_f(z)\,-\,2\overline{z}\,-\,(\alpha + 1)
\frac{1-\overline{z}}{1-z}\right|\,\leq\,\alpha - 1.
\ee
This proves the first part of the assertion in the theorem.
The second part follows from the fact that $|\omega(z)| = 1$ if and only if
$\omega(z)\equiv e^{i\theta}$, $\theta \in [0,2\pi],$ and that the solution of the
differential equation (\ref{p8eq1}) in this case is given by $f(z)=g_{\theta}(z).$
The relation between boundary points of the above circle and the extremal function
becomes clear from the identity
$$
(1-|z|^2)\frac{g_{\theta}''(z)}{g_{\theta}'(z)}\,-\,2\overline{z}\,-\,(\alpha + 1)
\frac{1-\overline{z}}{1-z}\,=\,(\alpha-1)\frac{e^{i\theta}+\overline{z}}{1+e^{i\theta}z}.
$$
This completes the proof of the theorem.
\epf

\br
We  remark here that for $f\in Co(\alpha)$, the sharp inequality
(\ref{p8eq5aa}) was obtained by Cruz and Pommerenke in \cite[ Theorem 3]{Pom-Cruz}.
Their result proves only a one way implication, namely
the condition on the disk of variability
of the pre-Schwarzian is necessary for $f$ to belong to $Co(\alpha)$.
In our theorem we have actually shown that this condition
is not only necessary for $f$ to belong to $Co(\alpha)$ but is also sufficient.
\er

\bcor
Let $f\in Co(\alpha), \alpha\in [1,2]$. Then,
$4\,\leq\,\|T_f\|\,\leq 2\alpha+2.$ The equality holds in lower estimate for the function
$g_\pi$ and in upper estimate for the function $g_0$ which are described in the statement
of the above theorem.
\ecor
\bpf
Since
$$\sup_{|z|=r}\left|\,2\overline{z}\,+\,(\alpha + 1)
\frac{1-\overline{z}}{1-z}\right|\,=\,2r\,+\,1\,+\,\alpha,
$$
where the maximum is attained for $z=r$, we deduce immediately from (\ref{p8eq5aa}),
that
$$
2+2r\,\leq \sup_{|z|\leq r}(1-|z|^2) |T_f(z)|\,\leq 2\alpha+2r.
$$
The lower bound is attained for $f=g_{\pi}$, and the upper bound for $f=g_0$. Indeed,
we see that
$$ \sup_{|z|\leq r}(1-|z|^2)|T_{g_{\pi}(z)}|
= 2 \sup_{|z|\leq r}\frac{1-|z|^2}{|1-z|}= 2(1+r)
$$
and
$$  \sup_{|z|\leq r}(1-|z|^2)|T_{g_{0}(z)}|
=  \sup_{|z|\leq r}\frac{|2\alpha+2z|(1-|z|^2)}{|1-z^2|}= 2(\alpha+r).
$$
Now,  letting $r \ra 1$, we get the  sharp estimates
$$4\,\leq\,\|T_f\|\,\leq 2\alpha+2,~~ f\in Co(\alpha).
$$
\epf

\br
It is well-known that for the class $\mathcal{K}$ of convex univalent functions $f$,
the pre-Schwarzian norm $\|T_f\|$ satisfies the sharp inequality
$\|T_f\|\leq 4$
and the equality holds for the convex function $g_{\pi}(z)=z/(1-z)$. Moreover,
we observe that $\|T_f\|\geq 4$ for the class of concave functions
and the equality holds for the function $g_{\pi}(z)=z/(1-z)$ which
is common to both the classes and the only function in $Co(\alpha)$ with $\alpha=1.$
\er

As a consequence of Theorem \ref{p8thm1}, we can obtain a distortion theorem.

\bthm[{Distortion Theorem}]\label{p8thm1a}
Let $\alpha\in (1,2]$. Then, for each $f\in Co(\alpha)$, we have
$$ \frac{(1-r)^{\alpha-1}}{(1+r)^{\alpha+1}}\leq |f'(z)|\leq \frac{(1+r)^{\alpha-1}}
{(1-r)^{\alpha+1}},~~ |z|=r<1.
$$
For each $z\in\ID, z\neq 0$, equality occurs if and only if $f=g_\theta$, where
$\theta\in [0, 2\pi)\setminus \{\pi\}$.
\ethm
\bpf
In view of the inequality  (\ref{p8eq5aa}), it follows easily that
$$
\left|\frac{zf''(z)}{f'(z)}-\frac{2r^2}{1-r^2}\right|\leq \frac{2\alpha r}{1-r^2},
\quad  \mbox{$|z|=r<1$}.
$$
A standard argument (see for eg. \cite[Theorem 2.5]{Duren}) gives the desired
estimate for $|f'(z)|$. Also the sharpness part is easy to verify and so, we skip
the routine calculation.
\epf

In order to include an inclusion result, we need to introduce
another notation. Let $H^{p}$,  $p\in (0, \infty)$, denote
the standard Hardy space of analytic functions on the unit disk $\ID$
(see for eg. Duren \cite[p.~60--62]{Duren}). It is wellknown that
$\mathcal{S}$ is included in $H^{p}$ for $0<p<1/2$. For the class of
convex functions, the range for $p$ can be extended to $0<p<1$.

\bcor
$Co(\alpha)\subset H^{p}$ for $0<p<1/\alpha$. The result is best possible.
\ecor
\bpf
We fix $z=re^{i\theta}$ with $0<r<1$. As $f(0)=0$, we observe
$$
f(z)=\int_{0}^{r}f'(\rho e^{i\theta}) e^{i\theta} d\rho.
$$
Hence by the distortion theorem and a mild computation, one has
$$
|f(z)|\leq \int_{0}^{r} \frac{(1+\rho)^{\alpha-1}}
{(1-\rho)^{\alpha+1}} d\rho \leq \frac{K}{(1-r)^{\alpha}},
$$
for some positive constant $K$. The desired result
follows form the last inequality and the Prawitz' theorem
(see for eg. \cite[Theorem 2.22]{Duren}).
\epf

There has been a number of investigations on basic subclasses of univalent
functions by fixing the second coefficient of functions in these classes.
Therefore, it is natural to obtain an analog of Theorem \ref{p8thm1}
for functions in $ f\in Co(\alpha)$ with fixed second coefficient.
Our next result gives the set of
variability of the functional $(1-|z|^2)T_f(z)$ for $f\in Co(\alpha)$
whenever $f''(0)$ is fixed.

\bthm \label{p8thm2}
Let $f\in Co(\alpha)$, $\alpha\in (1,2]$.
Then the set of variability of the functional $T_f(z)(1-|z|^2)$,  $f\in Co(\alpha)$,
whenever $f''(0)=\alpha + 1 +(\alpha - 1)a$ with $a \in \overline \ID$ being fixed,
is the disk
\beq%\label{p8eq10}
\nonumber
&&\left|(1-|z|^2)T_f(z)\,-\,2\overline{z}\,-\,(\alpha + 1)\frac{1-\overline{z}}{1-z}-
(\alpha-1)\frac{\overline{z}(1+|a|^2+\overline{a}\overline{z})+a}
{1+|z|^2+2 {\rm Re}\,(az)}\right|\\\nonumber
&& \leq (\alpha-1)\frac{(1-|a|^2)|z|}{1+|z|^2+2 {\rm Re}\,(az)}.\nonumber
\eeq
\ethm
\bpf
As in the proof of Theorem \ref{p8thm1}, a calculation reveals that for $f\in Co(\alpha)$,
\be\label{p8eq6}
(1-|z|^2)T_f(z)\,-\,2\overline{z}\,-\,(\alpha + 1)
\frac{1-\overline{z}}{1-z}=(\alpha-1)\frac{\overline z+ \omega(z)}{1+z\omega(z)},
\ee
where $\omega :\ID\to \overline{\ID}$ is an unimodular bounded analytic function.

We see that by (\ref{p8eq1}) fixing $f''(0)$
is equivalent to fixing $\omega(0)$, where $\omega$ is as above. Indeed we have
$f''(0)=\alpha + 1 +(\alpha - 1)\omega(0)$.
Now, let
\be\label{p8eq7}
\omega(z)\,=\,\frac{\omega(0)+z\phi(z)}{1+\overline{\omega(0)}z\phi(z)},
\ee
where $\phi:\ID\to \overline{\ID}$ is again an analytic
unimodular bounded function.
For convenience, we let $\omega(0)= a$. Then from (\ref{p8eq7}), we get
$$
z\phi(z)=\frac{\omega(z)-a}{1-\overline{a}\omega(z)}
$$
and a computation shows that $|\phi(z)|\leq 1$ if and only if
\be\label{p8eq8}
|\omega(z) - W_0|\leq R, ~z\in \ID,
\ee
where
$$
W_0= \frac{a(1-|z^2|)}{1-|z|^2|a|^2}\quad\mbox{and} \quad
R=  \frac{|z|(1-|a|^2)}{1-|z|^2|a|^2}.
$$
In order to complete the proof, we let
$$
W= \frac{\overline z+ \omega(z)}{1+z\omega(z)}.
$$
This gives
$$
\omega(z)= \frac{W-\overline z}{1-z W}
$$
so that (\ref{p8eq8}) is equivalent to
$$
\left|\frac{W-\overline z}{1-z W}- W_0\right|\leq R.
$$
By a routine calculation the last inequality reduces to
\be\label{p8eq9}
\left|W-\frac{(1+\overline{W_0}\overline{z})(\overline{z}+W_0)-R^2{\overline z}^2}
{|1+W_0 z|^2- R^2|z|^2}
\right|\leq \frac{R(1-|z|^2)}{|1+W_0z|^2- R^2|z|^2}.
\ee
An easy exercise gives
$$|1+W_0 z|^2- R^2|z|^2=
 \left(\frac{1-|z|^2}{1-|a|^2|z|^2}\right)(1+|z|^2+2 {\rm Re}\,(az)),
$$
and
$$(1+\overline{W_0}\overline{z})(\overline{z}+W_0)-R^2{\overline z}^2
= \frac{1-|z|^2}{1-|a|^2|z|^2}\left(\overline{z}(1+|a|^2+\overline{a}\overline{z})+a\right).
$$
Using the above two equalities, we see that the inequality (\ref{p8eq9}) takes the
following equivalent form
$$
\left|W-\frac{\overline{z}(1+|a|^2+\overline{a}\overline{z})+a}{1+|z|^2+2 {\rm Re}\,(az)}
\right|\leq \frac{(1-|a|^2)|z|}{1+|z|^2+2 {\rm Re}\,(az)}.
$$
Hence from (\ref{p8eq6}) we get that the set of variability of the functional
$T_f(z)(1-|z|^2)$ is
\beq%\label{p8eq10}
&&\left|(1-|z|^2)T_f(z)\,-\,2\overline{z}\,-\,(\alpha + 1)\frac{1-\overline{z}}{1-z}-
(\alpha-1)\frac{\overline{z}(1+|a|^2+\overline{a}\overline{z})+a}
{1+|z|^2+2 {\rm Re}\,(az)}\right|\\\nonumber
&& \leq (\alpha-1)\frac{(1-|a|^2)|z|}{1+|z|^2+2 {\rm Re}\,(az)}\nonumber
\eeq
(where $a=\omega(0)$ is fixed).
Whenever $f''(0)=\alpha + 1 +(\alpha - 1)e^{i\theta}$, i.e. $a= e^{i\theta}$, the last
inequality reduces to
$$
(1-|z|^2)T_f(z)\,-\,2\overline{z}\,-\,(\alpha + 1)\frac{1-\overline{z}}{1-z}=
(\alpha-1)e^{i\theta}\left(\frac{e^{i\theta}+\overline z}{1+ze^{i\theta}}\right)^2.
$$
As
$$
(1-|z|^2)\frac{g_{\theta}''(z)}{g_{\theta}'(z)}\,-\,2\overline{z}\,-\,(\alpha + 1)
\frac{1-\overline{z}}{1-z}\,=\,(\alpha-1)\frac{e^{i\theta}+\overline{z}}{1+e^{i\theta}z},
$$
the boundary of the disk of variability is attained if and only if
$f=g_\theta$ where $g_\theta$ is given in Theorem \ref{p8thm1}.
\epf

\bcor
Let $f\in Co(\alpha)$ and $f''(0)=\alpha+1$ be fixed. Then,
$$
3+\alpha\leq \|T_f\|\leq 2+2\alpha.
$$
\ecor
\bpf
Setting $a=\omega(0)=0$ in Theorem \ref{p8thm2} we get
$$
\left|(1-|z|^2)T_f(z)\,-\,2\overline{z}\,-\,(\alpha + 1)\frac{1-\overline{z}}{1-z}-
(\alpha-1)\frac{\overline{z}}{1+|z|^2}\right| \leq (\alpha-1)\frac{|z|}{1+|z|^2}.
$$
This inequality  easily  gives the required estimates for the pre-Schwarzian norm.
\epf

\section{Convolution Characterization and Coefficient Estimates}

If $f,g\in \mathcal{H}$, with
$$f(z)=\sum_{n=0}^\infty a_nz^n ~~\mbox{ and }~~
g(z)=\sum_{n=0}^\infty b_nz^n,
$$
then the Hadamard product (or convolution) of $f$ and $g$
is defined by the function
$$(f\star g)(z)=\sum_{n=0}^\infty a_nb_n z^n.
$$
Clearly, $f\star g\in \mathcal{H}$. In view of the Hadamard convolution, it is now
possible to present a new characterization for functions in the class $Co(\alpha)$.
The following result will be useful although we did not gain much inroads in this direction.

%We apply similar method as given in \cite{Sil-78} to derive
%the following theorem.
\bthm\label{p8thm3}
Let $1<\alpha\leq 2$. Then, $f\in Co(\alpha)$  if and only if
\be\label{p8eq11}
\frac{1}{z}\left[f \star \frac{(\alpha-1)z-(\alpha+1+2x)z^2}{(1-z)^3}\right]
+\left[f\star  \frac{((\alpha+1)x+2)z-(\alpha -1)xz^2}{(1-z)^3}\right]\neq 0
\ee
for all $|z|<1$ and for all $x$ with $|x|=1$. Equivalently, this holds if and only if
\be\label{p8eq12}
\sum_{n\geq 0}A_nz^n \neq 0, \quad A_0= \alpha-1,\, (z\in \ID,~ |x|=1)
\ee
where
\be\label{p8eq15a}
f(z)=z+\sum_{n\geq 2}a_nz^n
\ee
and
$$A_n= (\alpha-n-1-nx)(n+1)a_{n+1}+[n+1+(n+\alpha)x]na_n \quad  (n\geq 1,~ a_1=1).
$$
\ethm
\bpf
We recall
$f\in Co(\alpha)$ if and only if ${\rm Re}\,P_f(z)>0$  in $\ID$, where
$$P_f(z) =\frac{2}{\alpha-1} \left(\frac{\alpha+1}{2}
\frac{1+z}{1-z}-\frac{zg'(z)}{g(z)}\right)
$$
with  $g(z)=zf'(z)$. We note that $P_f$ is analytic in $\ID$ with $P_f(0)=1$.
Thus, $f\in Co(\alpha)$ is equivalent
to
$$ P_f(z) \neq \frac {x-1}{x+1}\quad (z\in \ID,~ |x|=1, ~x\neq -1)
$$
which, by a simplification, is  same as writing
%$$[\alpha+1-\frac{x-1}{x+1}(\alpha-1)]\frac{g(z)}{z} +
%[\alpha+1 +\frac{x-1}{x+1}(\alpha-1)]g(z)-2g'(z)+2zg'(z)\neq 0
%$$
%which is equivalent to
\beq\label{p8eq13}
\left (\frac{x+\alpha }{x+1}\right )\frac{g(z)}{z} +
\left (\frac{\alpha x +1}{x+1}\right )g(z)-g'(z)+zg'(z)\neq 0.
\eeq
Recall that
$$ \frac{g(z)}{z}= \frac{g(z)}{z} \star \frac{1}{1-z},  ~zg'(z)= g(z) \star
\frac{z}{(1-z)^2},~
f'(z) \star p(z)=\frac{f(z)}{z} \star (zp)'(z).
$$
Using these identities, $g(z)=zf'(z)$, (\ref{p8eq13}) gives that $f\in Co(\alpha)$  if
and only if
$$ f'(z) \star \left (\frac{(x+\alpha)(1-z)-(x+1)}{(1-z)^2}\right)
+z\left [f'(z)\star \left(\frac{(\alpha x+1)(1-z) +(x+1)}{(1-z)^2}\right)\right ]\neq 0.
$$
After some simplification the above takes the following equivalent form
\beq\label{p8eq14}\nonumber
&&\frac{f(z)}{z} \star \left [\frac{\alpha-1-(\alpha+1+2x)z}{(1-z)^3}\right ]\\
&&+z\left[\frac{f(z)}{z}\star \left (\frac{(\alpha+1)x+2-(\alpha-1)xz}{(1-z)^3}\right )
 \right]\neq 0
\eeq
which gives (\ref{p8eq11}). To obtain the series formulation of it, we first observe that
\beqq
p(z)&=&\frac{\alpha-1-(\alpha+1+2x)z}{(1-z)^3}\\
&=& \alpha-1+ \sum_{n\geq 1}(n+1) \left[\alpha-n-1-nx\right]z^n
\eeqq
and
\beqq
q(z)&=& \frac{((\alpha+1)x+2)z-(\alpha -1)xz^2}{(1-z)^3}\\
&=& ((\alpha+1)x+2)z+ \sum_{n\geq 2}n[n+1+(n+\alpha)x]z^n.
\eeqq
Using these two identities, (\ref{p8eq14}) can be written in terms of convolution as
follows:
$f\in Co(\alpha)$ $(1<\alpha\leq 2)$ if and only if
$$\frac{f(z)}{z} \star p(z) +f(z)\star  q(z) \neq 0
$$
which is same as (\ref{p8eq12}). We complete the proof.
\epf

In order to reveal the interaction between
the class $Co(\alpha)$ and wellknown Kaplan class, we need to
introduce the following definition.

\begin{defn}
A nonvanishing analytic function $s$ in $\ID$ with $s(0)=1$ is said to belong to the
{\it Kaplan class} $K(\alpha, \beta)$ ($\alpha\geq 0, \beta\geq 0$) if for $0<r<1$ and
$\theta_1<\theta_2<\theta_1+2\pi$ we have
$$
-\alpha\pi\leq \int_{\theta_1}^{\theta_2}\left\{{\rm Re}\,
\frac{re^{i\theta}s'(re^{i\theta})}{s(re^{i\theta})}-\frac{1}{2}(\alpha-\beta)\right\}
\,d\theta
\leq \beta\pi.
$$
\end{defn}
Following the notation of Sheil-Small \cite{Sheil-83}, for $\lambda$ real we consider
$$\Pi_\lambda =
\left\{\pay
\ba{rl}
K(\lambda, 0) & (\lambda \geq 0) \\
K(0, -\lambda)& (\lambda < 0).
\ea\right.
$$
This gives  $f\in \Pi_\lambda$ if and only if for $z\in \ID$,
$${\rm Re}\,
z\frac{s'(z)}{s(z)}
\left\{\pay
\ba{rl}
<\frac{1}{2}\lambda & (\lambda > 0) \\
>\frac{1}{2}\lambda & (\lambda < 0).\\
\ea\right.
$$
The class $\Pi_0=K(0,0)$ contains the constant function $s(z)=1$ only
(compare \cite{Sheil-83}).

%Following the same notation we get a nice connection
%of the functions in the class $Co(\alpha)$ and the Kaplan class which we describe below
\bthm\label{p8thm4}
Let $\alpha\in (1, 2]$.
A function $f\in Co(\alpha)$ if and only if there exists a function $s\in \Pi_{\alpha-1}$
such that
\be\label{p8eq15}
f(z)=\int_{0}^{z}\frac{s(t)}{(1-t)^{\alpha+1}}\, dt.
\ee
\ethm
\bpf
The function $f$ given by (\ref{p8eq15}) satisfies
$$
s(z)=(1-z)^{\alpha+1}f'(z).
$$
A computation from (\ref{p8eq1}) shows that
\be\label{p8eq16}
\frac{s'(z)}{s(z)}= \frac{f''(z)}{f'(z)}-\frac{\alpha+1}{1-z}=
\frac{\alpha-1}{2}\,\frac{1-P_f(z)}{z}.
\ee
This gives
$${\rm Re}\,\left(z\frac{s'(z)}{s(z)}\right)<\frac{\alpha-1}{2},~ z\in \ID
$$
if and only if
$$
{\rm Re}\,\left[\frac{(\alpha+1)}{2}\frac{1+z}{1-z}-1-z \frac{f''(z)}{f'(z)}\right]>0,
~ z\in \ID.
$$
The proof follows.
\epf

The functions $s$ defined as above produces a simple characterization
in terms of the Hadamard product.

\bthm
A function $f\in Co(\alpha)$ if and only if
$$ s(z)\star \left(\frac{z}{(1-z)^2}+\frac{1-\alpha}{x+1}\frac{1}{1-z}\right)
\neq 0,\quad z\in \ID,~ |x|=1, x\neq -1,
$$
for some $s\in \Pi_{\alpha-1}$.
\ethm
\bpf
We recall that, $f\in Co(\alpha)$ is equivalent
to
$$ P_f(z) \neq \frac {x-1}{x+1}\quad (z\in \ID,~ |x|=1, ~x\neq -1),
$$
which by (\ref{p8eq16}) is equivalent to
$$
\frac{-2}{\alpha-1}\left(\frac{zs'(z)}{s(z)}\right)+1 \neq \frac {x-1}{x+1}.
$$
A simplification gives
$$
(1-\alpha) \left(s(z)\star\frac{1}{1-z}\right)+ (x+1)\left(s(z)\star \frac{z}{(1-z)^2}\right)\neq 0
$$
and the desired condition follows.
\epf

Moreover, if $f\in Co(\alpha)$, we may define a function $\varphi: \ID\to \overline{\ID}$ by
\be\label{adde0}
P_f(z)=\frac{1+z\varphi(z)}{1-z\varphi(z)}.
\ee
In view of this and (\ref{p8eq16}), we have
\be\label{adde1}
s'(z)=\varphi(z)((1-\alpha)s(z)\,+\,zs'(z)).
\ee
We want to find bounds for the moduli of the Taylor coefficients $b_k$,
$k\in \IN,$ of the function $s$ that are defined via the series representation
$$s(z)=\sum_{k=0}^{\infty}b_kz^k,\quad b_0=1.
$$
To that end we use Theorem 2.2 in \cite{P} (compare also \cite{Cl,Ro, Rog}).
Using these methods, we see that (\ref{adde1}) implies the inequalities
$$\sum_{k=1}^N k^2|b_k|^2\leq
\sum_{k=0}^{N-1}(k+1-\alpha)^2|b_k|^2,\quad N\in \IN.
$$
Since $k+1-\alpha < k$, we may use mathematical induction to prove the inequality
\be\label{adde2}
|b_k|\leq\frac{\alpha-1}{k},\quad k\in \IN.
\ee
Equality in (\ref{adde2}) can be achieved if and only if
$$b_m=0,\quad m=1, \ldots ,k-1.
$$
If we insert this into  (\ref{adde1}) and assume equality in (\ref{adde2}), we
recognize that this implies
$$\varphi(z)=e^{i\theta}z^{k-1},\quad \theta\in[0,2\pi].
$$
Solving (\ref{p8eq1}) with $P_f$ defined by (\ref{adde0}) and this function
$\varphi$ leaves us with the fact that the functions
\be\label{adde3}
f'(z)=\frac{(1-e^{i\theta} z^k)^{\frac{\alpha-1}{k}}}{(1-z)^{\alpha+1}}
\ee
are the unique extremal functions for the inequalities (\ref{adde2}).
The Schwarz-Christoffel formula implies that the functions (\ref{adde3}) deliver
functions in $Co(\alpha)$ that map the unit disk conformally onto the complement
of an unbounded polygon with $k$ or $k-1$ finite vertices.

\br
Since we may consider the functional that maps $f$ into  $s^{(k)}(0)/k!$ as a
linear functional on  $Co(\alpha)$ with the set of variability described
by (\ref{adde2}) with a unique extremal function corresponding to any boundary
point, we get new examples supporting a conjecture formulated in \cite{BPW}.
There, we conjectured that any conformal map of the unit disk onto  the
complement of an unbounded convex polygon is an extremal point of the closed
convex hull of $Co(\alpha)$.
\er

In view of the discussion above  one can quickly get the following
\bthm\label{p8thm5}
Let $f\in Co(\alpha)$ have the expansion $(\ref{p8eq15a})$. Then the following sharp
inequality holds
\be\label{p8eq10a}
\left|\sum_{k=0}^{n}(-1)^{k} {\alpha+1\choose n-k}(k+1)a_{k+1} \right|
\leq \frac{\alpha-1}{n} \quad (n\geq 1).
\ee
In particular, we have
\bee
\item[(i)]
$\ds
\left|a_2- \frac{\alpha+1}{2}\right|\leq \frac{\alpha-1}{2}
$
\item[(ii)]
$\ds
\left|3a_3-2(\alpha+1)a_2+\frac{\alpha(\alpha+1)}{2}\right|\leq \frac{\alpha-1}{2}.
$
\eee
\ethm
\bpf
We deduce from Theorem \ref{p8thm4} that $f\in Co(\alpha)$ if and only if $s(z)\in \Pi_{\alpha-1}$
where
$$s(z)=(1-z)^{\alpha+1}f'(z).
$$
Comparing the coefficient $z^n$ in the series expansion of the functions in involved above,
we obtain that
$$b_n=\sum_{k=0}^{n}(-1)^{k} {\alpha+1\choose n-k}(k+1)a_{k+1}
$$
and the desired inequality (\ref{p8eq10a}) follows, if we use the
estimate (\ref{adde2}). The two particular cases follow,
if we let $n=1,2$. Also the estimate is sharp
for each $n$ for the functions $f\in Co(\alpha)$ such that
$$
f'(z)=\frac{(1-e^{i\theta} z^n)^{\frac{\alpha-1}{n}}}{(1-z)^{\alpha+1}}.
$$
%The equality in the particular cases follows for functions defined above
%if we choose $n=1, 2$ for the
%estimates (i) and (ii) respectively.
\epf

\br
Case (i) of Theorem \ref{p8thm5} is wellknown whereas Case (ii) of
Theorem \ref{p8thm5} gives that for $f\in Co(2)$ one has
$$\left|1-2a_2+a_3\right|\leq \frac{1}{6}.
$$
This result is obtained recently in \cite[Theorem 3]{Wir-05a}.
\er

The classical Alexander transform $\int_0^z(f(t)/t)\,dt$ provides a one-to-one
correspondence between $\mathcal{S}^*$ and $\mathcal{K}$. It is then natural to
ask whether a similar correspondence can be established for $\mathcal{S}^*$
and $Co(\alpha)$.  The answer is provided by the following result
which may be used to study the geometric properties of $\Lambda _\phi (z)$
when $\phi $ belongs to various subclasses of $\mathcal{S}$.

%result which is indeed an equivalent form of Theorem \ref{p8thm4}.
\bthm \label{thm1}
Let $\alpha\in (1, 2]$. A function $f\in Co(\alpha)$ if and only
if there exists a $\phi\in \mathcal{S}^{*}$ such that $f(z)=\Lambda _\phi (z)$,
where
$$\Lambda _\phi (z)
=\int_{0}^{z}\frac{1}{(1-t)^{\alpha+1}}\left(\frac{t}{\phi(t)}\right)^{(\alpha-1)/2}\, dt.
$$
\ethm\bpf In view of Theorem~\ref{p8thm4}, it suffices to  show that
$s\in\Pi_{\alpha-1}$ if and only if there exists a $\phi\in \mathcal{S^{*}}$ such
that
$$ \phi (z) =z\left(s(z)\right)^{2/(1-\alpha)}, ~~ z\in \ID.
$$
However, this fact is clear because
$$\frac{2}{\alpha-1} {\rm Re}\,\left(z\frac{s'(z)}{s(z)}\right)
= 1- {\rm Re}\,\left(z\frac{\phi'(z)} {\phi(z)}\right),~~ z\in \ID.
$$
\epf

For $f\in Co(\alpha)$, the above characterization shows that there exists $\phi\in
\mathcal{S}^{*}$ such that
\be\label{p8eq15b}
f'(z)=\frac{1}{(1-z)^{\alpha+1}}\left(\frac{z}{\phi(z)}\right)^{(\alpha-1)/2}.
\ee
Set
$$f(z)=z+\sum_{n\geq 2}a_nz^n ~\mbox{ and }~ \phi(z)=z+\sum_{n\geq 2}\phi_nz^n.
$$
A comparison of coefficients of $z^2$ on both side of (\ref{p8eq15b}) yields
$$ 3a_3=-\frac{\alpha-1}{2}\phi _3+\frac{(\alpha+1)(\alpha-1)}{8}\phi_2^{2}
-\frac{\alpha^2-1}{2}\phi_2 +\frac{(\alpha+1)(\alpha+2)}{2}.
$$
That is,
$$ \frac{3}{\alpha^2-1}\left(a_3-\frac{(\alpha+1)(\alpha+2)}{6}\right)
=A(\phi_3,\phi_2,\alpha)
$$
where
$$A(\phi_3,\phi_2,\alpha)=-\frac{1}{2(\alpha+1)}\left[\phi_3-\frac{\alpha+1}{4}\phi_2^2+(\alpha+1)\phi_2\right].
$$
In view of the result of Avkhadiev and Wirths \cite[Corollary 3]{Avk-Wir-05}, one
has the following result which does not seem to be known in this form.

\bcor
Let $\phi(z)=z+\sum_{n\geq 2}\phi_nz^n$ belong to $\mathcal{S}^{*}$ and $\alpha\in (1, 2]$.
Then $A(\phi_3,\phi_2,\alpha )\in \overline{h(\ID)}$, where
$$h(z)=z+\frac{\alpha-2}{2(\alpha+1)}z^2.
$$
\ecor

%It would be interesting to generalize this result (compare with the
%classical Fekete-Szeg\"{o} type problem) and derive analog of it for the
%class of convex or the class of close-to-convex functions is unknown.

It is not clear whether the present restriction on $\alpha$ is essential
in the last corollary.


\begin{thebibliography}{99}
\bibitem{Avk-Wir-06} {\sc F.G. Avkhadiev, Ch. Pommerenke, and K.-J. Wirths}:
Sharp inequalities for the coefficient of concave schlicht functions,
\textit{Comment. Math. Helv.} {\bf 81}(2006), 801--807.

\bibitem{Avk-Wir-05} {\sc F.G. Avkhadiev  and K.-J. Wirths}:
Concave schlicht functions with bounded opening angle at infinity,
\textit{ Lobachevskii J. Math.} {\bf 17}(2005), 3--10.

%\bibitem{AW-07} {\sc F.G. Avkhadiev and K.-J. Wirths}:
%A proof of Livingston conjecture,
%\textit{Forum Math.} {\bf 19}(2007),  149--158.

\bibitem{Bek-72} {\sc J. Becker}:
L\"ownersche Differentialgleichung und quasikonform fortsetzbare schlichte
Funktionen (German),
\textit{J. Reine Angew. Math.} {\bf 255}(1972), 23--43.

\bibitem{BPW} {\sc B. Bhowmik, S. Ponnusamy, and K.-J. Wirths}:
Unbounded convex polygons, Blaschke products and concave schlicht functions,
\textit{Indian J. Math.} {\bf 50}(2008), 339--349.

%\bibitem{Choi-05} {\sc J. H. Choi, Y.C. Kim, S. Ponnusamy, T. Sugawa}:
%Norm estimates for Alexander transforms of convex functions of order $\alpha$,
%\textit{J. Math. Anal. Appl.} {\bf 303}(2005),  661--668.

\bibitem{Cl} {\sc J. G. Clunie:}  On meromorphic schlicht functions,
\textit{J. London Math. Soc.} {\bf 34}(1958), 215--216.

\bibitem{Pom-Cruz} {\sc L. Cruz  and Ch. Pommerenke}:
On concave univalent functions,
\textit{Complex Var. Elliptic Equ.} {\bf 52}(2007),  153--159.

\bibitem{Duren} {\sc P.L. Duren}:
Univalent functions, Springer-Verlag, 1983.


%\bibitem{Goodman} {\sc A.W. Goodman}:
%Univalent functions, vols. I and II,
%\textit{Mariner Publishing Co.,} Tampa, Florida, 1983.

%\bibitem{Hallen-Living-76} {\sc D.J. Hallenbeck and  A.E. Livingston}:
%Applications of extreme point theory to classes of multivalent functions,
%\textit{Trans. Amer. Math. Soc.} {\bf 221}(1976),  339--359.

\bibitem{H} {\sc H. Hornich}:
Ein Banachraum analytischer Funktionen in Zusammenhang mit den schlichten Funktionen,
\textit{Monatsh. Math.} {\bf 73}(1969), 36--45.

%\bibitem{Kustner-02} {\sc R. K\"{u}stner}:
%Mapping properties of hypergeometric functions and convolutions of starlike or
%convex functions of order $\alpha$,
%\textit{Comput. Methods Funct. Theory} {\bf 2}(2002),  597--610.

\bibitem{Kim-2002} {\sc Y. C. Kim and T. Sugawa}:
Growth and coefficient estimates for uniformly locally univalent functions
on the unit disk,
\textit{Rocky Mt. J. Math.} {\bf 32}(2002), 179--200.

\bibitem{Mil-Mocanu} {\sc S. S. Miller and P. T. Mocanu}:
Differential subordinations,
Marcel Dekker Inc. New York, Basel, 2000.

%\bibitem{Nez-05} {\sc  I. R. Nezhmetdinov, S. Ponnusamy}:
%New coefficient conditions for the starlikeness
%of analytic functions and their applications,
%\textit{Houston J. Math.} {\bf 31}(2005),  587--604. {Duren,P}

\bibitem{P} {\sc Ch. Pommerenke}: Univalent functions, Vandenhoeck and
Ruprecht, G\"ottingen, 1975.

\bibitem{Ro} {\sc M. S. Robertson}: Quasi-subordination and coefficient conjectures,
\textit{Bull. Amer. Math. Soc.} {\bf 76}(1970), 1--9.
\bibitem{Rog} {\sc W. W. Rogosinski}: On the coefficients of subordinate functions,
\textit{Proc. London Math. Soc.} {\bf 48}(1943), 48--82.
%\bibitem{Sil-78} {\sc H. Silverman, E. M. Silvia, D. Telage  }:
%Convolution conditions for convexity starlikeness and spiral-likeness,
%\textit{Math. Z.} {\bf 162}(1978), 125--130.
\bibitem{Sheil-83} {\sc T. Sheil-Small}: Coefficients and integral means of some classes of
analyic functions, \textit{Proc. Amer. Math. Soc.} {\bf 88}(1983), 275--282.
\bibitem{Wir-05a} {\sc K.-J. Wirths}:
Julia's lemma and concave schlicht functions,
\textit{Quaest. Math.} {\bf 28}(2005),  95--103.
\end{thebibliography}
\end{document}